\documentclass [12pt]{amsart}

\newtheorem{theorem}{Theorem}[section]
\newtheorem{lemma}[theorem]{Lemma}
\newtheorem{prop}[theorem]{Proposition}

\theoremstyle{definition}

\theoremstyle{remark}

\numberwithin{equation}{section}

\begin{document}

\title{Holomorphic Sobolev spaces, Hermite and special Hermite Semigroups and a
Paley-Wiener theorem for the windowed Fourier transform}
%Information for second author
\author{R. Radha}
\address{Department of Mathematics, Indian Institute of Technology, Chennai-600 036, India.}
\email{radharam@iitm.ac.in}

%Information for first author
\author{S. Thangavelu}
%    Address of record for the research reported here
\address{Department of Mathematics, Indian Institute of Science, Bangalore-560012, India.}
 {\email{veluma@math.iisc.ernet.in}
%    \thanks will become a 1st page footnote.
%\thanks{The first author was supported in part by NSF Grant \#000000.}

%\thanks{Support information for the second author.}

%    General info

%\date{January 1, 2001 and, in revised form, June 22, 2001.}

%\dedicatory{This paper is dedicated to our advisors.}


\begin{abstract}
The images of Hermite and Laguerre Sobolev spaces under the Hermite
and special Hermite semigroups (respectively) are characterised.
These are used to characterise the image of Schwartz class of rapidly
decreasing functions $f$ on $\mathbb{R}^n$ and $\mathbb{C}^n$ under these semigroups. The image of the space of tempered distributions is also considered
and a Paley-Wiener theorem for windowed ( short-time) Fourier transform is
proved.
\end{abstract}

%\maketitle
\maketitle \pagestyle{myheadings} \markboth{R.Radha and
S.Thangavelu}{Holomorphic Sobolev spaces}

%\section*{This is an unnumbered first-level section head}
%This is an example of an unnumbered first-level heading.

%% The correct journal style for \specialsection is all uppercase; a known bug
%% in amsart.cls prevents this, so input must be uppercase until it is fixed.
%\specialsection*{This is a Special Section Head}
%\specialsection*{THIS IS A SPECIAL SECTION HEAD}
%This is an example of a special section head%
%%%%%%%%%%%%%%%%%%%%%%%%%%%%%%%%%%%%%%%%%%%%%%%%%%%%%%%%%%%%%%%%%%%%%%%%
%\footnote{Here is an example of a footnote. Notice that this
%footnote text is running on so that it can stand as an example of
%how a footnote with separate paragraphs should be written.
%\par
%And here is the beginning of the second paragraph.}%

%%%%%%%%%%%%%%%%%%%%%%%%%%%%%%%%%%%%%%%%%%%%%%%%%%%%%%%%%%%%%%%%%%%%%%%%
\vspace{-.5cm}

\section{Introduction}
A classical result of Bargmann and Fock states that the image of
$L^2(\mathbb{R}^n)$ under the Gauss-Weierstrass semigroup can be
described as a weighted Bergman space of entire functions. Such
results are known for the semigroups generated by the
Laplace-Beltrami operators on compact Lie groups and more generally
on compact symmetric spaces. (See \cite {Hal94} and \cite {Ste99}).
Similar results have been  proved in the literature for Hermite and
special Hermite semigroups as well. (See  \cite {Byu93}, \cite
{Kro05} and \cite {Tha06} for further details). However, unlike the
above cases, it was shown by Kr\"otz, Thangavelu and Xu in \cite
{Kro05} that in the case of Heisenberg group $\mathbb{H}^n,$ the
image of $L^2(\mathbb{H}^n)$ under the heat kernel transform is not
a weighted Bergman space but  the direct sum of two weighted Bergman
spaces. Bargmann in \cite {Bar61} obtained a characterization for
the image of Schwartz class of rapidly decreasing functions on
$\mathbb{R}^n$ under the Segal-Bargmann transform. He showed that if
$F$ is a holomorphic function on $\mathbb{C}^n,$ then there exists a
function $f\in \mathcal{S}(\mathbb{R}^n)$ with $F=C_t f$ ($C_t$
denotes the Segal-Bargmann transform) if and only if $F$ satisfies
\newpage
$$|F(x+iy)|^2\leq A_n \frac{e^{\frac{y^2}{t}}}{[1+(x^2+y^2)]^{2n}}$$
for some sequence of
constants $A_n, n=1,2,3,\cdots.$ Hall and Lewkeeratiyutkul in
\cite {Hal04} characterized the image of Sobolev spaces under the
Segal- Bargmann transform on a compact Lie group $K.$ They used
this result to obtain a characterization for the functions in
the image of $C^{\infty}(K)$ under this transform. Using Gutzmer's
formula the images of Sobolev spaces under the Segal-Bargmann
transform on compact Riemannian symmetric spaces were
characterised by Thangavelu recently in \cite{Tha07} extending
the results of \cite{Hal04}.

%In \cite {Tha07}, Thangavelu characterized the image of
%Sobolev spaces under the Segal-Bargmann transform on a compact
%Riemannian symmetric space $X.$ In that paper he also
%characterized the image of $C^{\infty}(X)$under such transform.
%In fact, this result extends the result of \cite {Hal04} to all compact
%symmetric spaces.

In this paper we characterize the image of
Hermite Sobolev spaces under the Hermite semigroup and Laguerre
Sobolev spaces under special Hermite semigroup. These results are then used
to  characterize  the images of Schwartz space on $\mathbb{R}^n$ and
$\mathbb{C}^n$ respectively under  these semigroups. Throughout this paper,
if $x =(x_1,x_2,\cdots,x_n)\in \mathbb{R}^n$ then $x^2$ will denote
$x_1^2+x_2^2+\cdots+x_n^2.$

\section{Holomorphic Hermite Sobolev spaces}
Let $h_k(x)=(2^k k!\pi^{\frac{1}{2}})^{-\frac{1}{2}} (-1)^k
\frac{d^k}{dx^k}(e^{-x^2})e^{\frac{x^2}{2}}, k=0,1,2,...$ denote
the normalised  Hermite functions. The multi-dimensional Hermite
functions are defined as follows: For $\alpha\in \mathbb{N}^n,
x\in \mathbb{R}^n,$ let $\Phi_{\alpha}(x)=\Pi_{i=1}^n
h_{\alpha_i}(x_i).$ Then the collection $\{\Phi_{\alpha}: \alpha
\in \mathbb{N}^n\}$ forms an orthonormal basis for
$L^2(\mathbb{R}^n).$ These functions $\Phi_{\alpha}~(\alpha\in
\mathbb{N}^n)$ are eigenfunctions of the Hermite operator $H=-
\Delta +\frac{1}{2}x^2.$ Any $f\in L^2(\mathbb{R}^n),$ has the
Hermite expansion $f=\sum_{\alpha}(f,\Phi_{\alpha})\Phi_{\alpha}.$
If $P_k$ denotes the orthogonal projection of
$L^2({\mathbb{R}^n})$ onto the eigenspace spanned by
$\{\Phi_{\alpha}:|\alpha|=k\},$ then the Hermite expansion can be
written as $f=\sum_{k=0}^{\infty}P_k f.$ For $f\in
L^2(\mathbb{R}^n),$ the series converges in the norm and for other classes of
functions $f,$ this denotes the formal sum.

The Hermite operator $H$ defines a semigroup, called the Hermite
semigroup denoted by $e^{-tH},t>0$ by the expansion
$$e^{-tH}f=\sum_{k=0}^{\infty}e^{-(2k+n)t}P_k f,
f\in L^2(\mathbb{R}^n).$$ On a dense subspace, we can write
$e^{-tH}$ as an integral operator with kernel $K_t(x,u),$
$$e^{-tH}f(x)=\int_{\mathbb{R}^n}f(u)K_t(x,u)du.$$
Using Mehler's formula (See eqn (1.1.36) in \cite {Tha93}), $K_t(x,u)$ can be explicitly written as
$$K_t(x,u)=(2\pi)^{-\frac{n}{2}}(\sinh(2t))^{-\frac{n}{2}}e^{-\frac{1}{2}coth(2t)(x^2+u^2)+\frac{1}{\sinh2t}x.u}.$$
The Hermite semigroup initially defined on $L^2\cap
L^p(\mathbb{R}^n)$ extends to the whole of $L^p(\mathbb{R}^n)$ and
$||e^{-tH}f||_p\leq c_t||f||_p, 1\leq p\leq {\infty}.$ It is easy
to see that $K_t(x,u)$ can be extended to $ \mathbb{C}^n $ as an entire
function
$K_t(z,w)$ and hence $e^{-tH}f$ can also be extended to $ \mathbb{C}^n $ and
the entire extension will be simply denoted by $e^{-tH}f(z), z=x+iy.$
It is well known that the image of
$L^2(\mathbb{R}^n)$ under the Hermite semigroup is a weighted Bergman space.
 To be more precise,
let $H_t(\mathbb{C}^n)$ denote the space of all entire functions
on $\mathbb{C}^n$ which are square integrable with respect to the
weight function
$U_t(x+iy)=2^n(\sinh4t)^{-\frac{n}{2}}e^{\tanh(2t)x^2-\coth(2t)y^2}.$
Then one has the following result which is due to Byun \cite{Byu93}.
\begin{theorem}
The Hermite semigroup $e^{-tH}: L^2(\mathbb{R}^n)\rightarrow
H_t(\mathbb{C}^n)$ is an isometric isomorphism.
\end{theorem}

The Hermite functions $\Phi_\alpha(x)
$ have extension to $\mathbb{C}^n$
as entire functions $\Phi_\alpha(z),$ called complexified Hermite functions.
They satisfy the orthogonality property
\begin{equation}
\int_{\mathbb{C}^n}\Phi_{\alpha}(z)\overline{\Phi_{\beta}(z)}U_t(z)dz=e^{2(2|\alpha|+n)t}\delta_{\alpha,\beta}.
\end{equation}
Let $\widetilde{\Phi}_{\alpha}(z)= e^{-(2|\alpha|+n)t}{\Phi_{\alpha}}(z).$ Then
$\{\widetilde{\Phi}_{\alpha}:\alpha \in \mathbb{N}^n\}$ forms an
orthonormal basis for $H_t(\mathbb{C}^n).$ Thus any $F \in
H_t(\mathbb{C}^n)$ can be written as $F = \sum_\alpha (F,
\tilde\Phi_\alpha)_{H_t(\mathbb{C}^n)} \tilde\Phi_\alpha.$

Hermite Sobolev spaces were studied by Thangavelu in \cite {Tha993}
in connection with regularity properties of twisted spherical means.
In order to give the definition, we consider the spectral
decomposition of $H:$
$$Hf=\sum_{k=0}^{\infty}(2k+n)P_k f.$$
Then the Hermite Sobolev space $W_{H}^{m,2}(\mathbb{R}^n)$ is
defined to be the image of $L^2(\mathbb{R}^n)$ under $H^{-m},$
where $m$ is a non negative integer. In other words, we say that
$f\in W_{H}^{m,2}(\mathbb{R}^n)$ if and only if
$$||f||_{W_{H}^{m,2}}^2=\sum_{k=0}^{\infty}(2k+n)^{2m}||P_k f||_2^2 <\infty.$$
The norm can be written as
$$||f||_{W_{H}^{m,2}}^2=\sum_{\alpha}(2|\alpha|+n)^{2m}|(f,\Phi_{\alpha})|^2.$$
The Sobolev space $W_{H}^{m,2}(\mathbb{R}^n)$ is a Hilbert space
under the inner product
$$(f,g)_{W_{H}^{m,2}(\mathbb{R}^n)}=\sum_{\alpha}(2|\alpha|+n)^{2m}(f,\Phi_{\alpha})\overline{(g,\Phi_{\alpha})}.$$
As $H^m f=\sum_{\alpha}(2|\alpha|+n)^m (f,\Phi_{\alpha})\Phi_{\alpha},$ the
above inner product can also be rewritten as
$$(f,g)_{W_{H}^{m,2}(\mathbb{R}^n)}=(H^m f,
H^m g)_{L^2(\mathbb{R}^n)}.$$

We shall now define the holomorphic Sobolev space
$W_t^{m,2}(\mathbb{C}^n)$ to be the image of
$W_{H}^{m,2}(\mathbb{R}^n)$ under $e^{-t H}.$ The space
$W_t^{m,2}(\mathbb{C}^n)$ is made into a Hilbert space simply
by transferring the Hilbert space structure of
$W_{H}^{m,2}(\mathbb{R}^n)$ to $W_t^{m,2}(\mathbb{C}^n)$ so that
the Hermite semigroup $e^{-t{H}}$ is an isometric isomorphism from
$W_{H}^{m,2}(\mathbb{R}^n)$ onto $W_t^{m,2}(\mathbb{C}^n).$
This means that
$$(F,G)_{W_t^{m,2}(\mathbb{C}^n)}=\sum_{\alpha}(2|\alpha|+n)^{2m}
(f,{\Phi}_{\alpha})_{L^2(\mathbb{R}^n)}
\overline{(g,{\Phi}_{\alpha})}_{L^2(\mathbb{R}^n)}$$
whenever $F=e^{-tH}f$ and $ G = e^{-t H}g $.
%If we denote $\widetilde{H}^m F=\sum_{\alpha}(2|\alpha|+n)^{2m}
%|(F,\widetilde{\Phi}_{\alpha})_{L^2(\mathbb{C}^n)}|^2$ then
%$$(F,G)_{W_t^{m,2}(\mathbb{C}^n)}=(\widetilde{H}^m
%F,\widetilde{H}^m G)_{L^2(\mathbb{C}^n)}.$$

Let $O(\mathbb{C}^n)$ denote the collection of all holomorphic
functions on $\mathbb{C}^n.$ Let $\mathcal{F}_t^m(\mathbb{C}^n)$
denote the space of all functions in $O(\mathbb{C}^n)$
which are square integrable with respect to the measure $|\frac
{d^{2m}}{dt^{2m}}U_t(z)| dz.$ We equip
$\mathcal{F}_t^m(\mathbb{C}^n)$ with the sesquilinear form
\begin{equation}
(F,G)_m=\int_{\mathbb{C}^n}F(z)\overline{G(z)}\frac
{d^{2m}}{dt^{2m}}U_t(z)dz.
\end{equation}
We shall show below that this defines a pre-Hilbert space structure
on $\mathcal{F}_t^m(\mathbb{C}^n).$ Let
$\mathcal{B}_t^m(\mathbb{C}^n)$ denote the completion of
$\mathcal{F}_t^m(\mathbb{C}^n)$ with respect to the norm induced by
the above inner product. In the following proposition, we also show
that $||F||_m$ and $||F||_{W_t^{m,2}}$ coincide up to a constant
multiple.
\begin{prop}
The sesquilinear form
$$(F,G)_m=\int_{\mathbb{C}^n}F(z)\overline{G(z)}\frac
{d^{2m}}{dt^{2m}}U_t(z)dz$$ is an inner product on
$\mathcal{F}_t^m(\mathbb{C}^n)$ and hence induces
a norm $||F||_m^2=(F,F)_m.$ We also have
$||F||_m^2=2^{2m}||F||_{W_t^{m,2}}$ for all functions $F = e^{-tH}f,$ with
$f \in \mathcal{S}(\mathbb{R}^n)$.
\end{prop}
\begin{proof}
Consider the integral
\begin{eqnarray}
\int_{\mathbb{C}^n}|F(x+iy)|^2 U_t(x+iy) dx dy.
\end{eqnarray}
Since the restriction of $F$ to $ \mathbb{R}^n $ has an orthogonal
expansion in terms of $\Phi_{\alpha},$ we have
$F(x+iy)=\sum_{\alpha}(F,\Phi_{\alpha})\Phi_{\alpha}(x+iy)$ where
$(F,\Phi_{\alpha})=\int_{\mathbb{R}^n}F(x)\Phi_{\alpha}(x)dx.$ Using
the orthogonality relation (2.1), we can show that the integral
(2.3) can be written as
$$\sum_{\alpha}|(F,\Phi_{\alpha})|^2 e^{2(2|\alpha|+n)t}.$$ By definition
$$ (F,F)_m =  \int_{\mathbb{C}^n}|F(x+iy)|^2\frac
{d^{2m}}{dt^{2m}}U_t(z)dz $$ which is the same as
$$ \frac {d^{2m}}{dt^{2m}}\sum_{\alpha}|(F,\Phi_{\alpha})|^2
e^{2(2|\alpha|+n)t}$$
$$ =  2^{2m}\sum_{\alpha}
(2|\alpha|+n)^{2m}|(F,\Phi_{\alpha})|^2 e^{2(2|\alpha|+n)t} $$ and
hence non-negative. Thus it follows that the sesquilinear form
defined in (2.2) is positive definite and induces the norm
$||F||_m.$ On the other hand if $F\in H_t{(\mathbb{C}^n)},$ we have
the expansion $$
F(z)=\sum_{\beta}(F,\widetilde{\Phi}_{\beta})_{H_t(\mathbb{C}^n)}\widetilde{\Phi}_{\beta}(z)$$
so that the restriction of $F$ to $\mathbb{R}^n$ can be written as
$$
F(x)=\sum_{\beta}(F,\widetilde{\Phi}_{\beta})_{H_t(\mathbb{C}^n)}\widetilde{\Phi}_{\beta}(x).$$
Thus
\begin{eqnarray*}
(F,\Phi_{\alpha})&=&
\int_{\mathbb{R}^n}\sum_{\beta}(F,\widetilde{\Phi}_{\beta})_{H_t(\mathbb{C}^n)}\widetilde{\Phi}_{\beta}(x)\Phi_{\alpha}(x)dx\\
&=&\int_{\mathbb{R}^n}\sum_{\beta}(f,\Phi_{\beta})e^{-(2|\beta|+n)t}\Phi_{\beta}(x)\Phi_{\alpha}(x)dx\\
&=&(f,\Phi_{\alpha})e^{-(2|\alpha|+n)t}
\end{eqnarray*}
where $F=e^{-tH}f.$
Again using (2.1), we get
\begin{eqnarray*}
||F||_m^2&=&2^{2m}\sum_{\alpha}
(2|\alpha|+n)^{2m}|(F,\Phi_{\alpha})|^2 e^{2(2|\alpha|+n)t}\\&=&
2^{2m}\sum_{\alpha} (2|\alpha|+n)^{2m}|(f,\Phi_{\alpha})|^2 \\&=&
2^{2m}\sum_{\alpha}
(2|\alpha|+n)^{2m}|(F,\widetilde\Phi_{\alpha})_{H_t(\mathbb{C}^n)}|^2\\&=&
2^{2m}||F||_{W_t^{m,2}}^2.
\end{eqnarray*}
\end{proof}

Using this proposition we can easily prove the following result on the image
of  Hermite Sobolev spaces under the Hermite semigroup.

\begin{theorem}
For every non negative integer m, $W_t^{m,2}(\mathbb{C}^n)$
coincides with $\mathcal{B}_t^m(\mathbb{C}^n)$ and the Hermite
semigroup $e^{-tH}$ is an isometric isomorphism of
$W_H^{m,2}(\mathbb{R}^n)$ onto $\mathcal{B}_t^m(\mathbb{C}^n)$ up
to a constant multiple.
\end{theorem}
\begin{proof} Since $\frac {d^{2m}}{dt^{2m}}U_t(z)$ is a polynomial 
in $ x,y $ times $U_t(z)$, any $ F $ which is square integrable with respect 
to $\frac {d^{2m}}{dt^{2m}}U_t(z)$ belongs to $\mathcal{H}_t(\mathbb{C}^n)$ 
and hence of the form $ e^{-t H}f .$ Further, it follows from the above 
proposition, as the norms $||F||_m$ and $||F||_{W_t^{m,2}}$
coincide, $ f \in W_H^{m,2}(\mathbb{R}^n).$ Consequently,
$\mathcal{F}_t^m(\mathbb{C}^n)$ is contained in
$W_t^{m,2}(\mathbb{C}^n).$ Notice that
$\widetilde{\Phi}_{\alpha}\in \mathcal{B}_t^m(\mathbb{C}^n).$
Further if $(F,\widetilde\Phi_{\alpha})_{W_t^{m,2}}$ $=0~\forall
~\alpha \in \mathbb{N}^n,$ then it can be easily seen that
$(F,\widetilde{\Phi}_{\alpha})_{H_t(\mathbb{C}^n)}=0 ~ \forall
~\alpha \in \mathbb{N}^n$ which forces
that $F=0.$ Hence $\mathcal{F}_t^m(\mathbb{C}^n)$ is dense in
$W_t^{m,2}(\mathbb{C}^n).$
\end{proof}

\section{Holomorphic Laguerre Sobolev spaces}

The special Hermite functions are defined using Hermite functions
by
$$\Phi_{\alpha\beta}(z)=(2\pi)^{-\frac{n}{2}}
\int_{\mathbb{R}^n}e^{ix\xi}\Phi_{\alpha}
(\xi+\frac{1}{2}y)\Phi_{\beta}(\xi-\frac{1}{2}y)
d\xi
$$
for $\alpha, \beta\in \mathbb{N}^n, z=x+iy.$ Then
$\{\Phi_{\alpha \beta}:\alpha, \beta\in \mathbb{N}^n\}$ forms an
orthonormal basis for $L^2(\mathbb{C}^n).$ The functions
$\Phi_{\alpha \beta}$ are eigenfunctions of the special Hermite
operator $L$ with eigenvalues $(2|\beta|+n),$ where
$$L=-\Delta_z+\frac{1}{4}|z|^2-i\sum_{i=1}^n \left( x_j\frac{\partial}{\partial y_j}-y_j\frac{\partial}{\partial
x_j}\right) $$where $-\Delta_z$ denotes the Laplacian on $\mathbb{C}^n.$
Any $f\in L^2({\mathbb{C}^n})$ has the special Hermite expansion
given by $f=\sum_{\alpha}\sum_{\beta}(f,\Phi_{\alpha
\beta})_{L^2(\mathbb{C}^n)}\Phi_{\alpha \beta}.$

Using the operation twisted convolution defined by
$$(f\times
g)(z)=\int_{\mathbb{C}^n}f(z-w)g(w)e^{\frac{i}{2}Im(z\bar w)}dw,$$
 the special Hermite expansion can be put in the following compact form:
$f=(2\pi)^{-n}\sum_{k=0}^{\infty}f\times \varphi_k$ where
$\varphi_k$ stands for the Laguerre function
$\varphi_k(z)=L_k^{n-1}(\frac{1}{2}|z|^2)e^{-\frac{1}{4}|z|^2}$ in
which $L_k^{n-1}$ denotes the $k$th Laguerre polynomial of type $n-1.$ For various results concerning Hermite and special Hermite expansions we refer to
Thangavelu \cite {Tha93}.

The special Hermite operator $L$ defines a semigroup ( called the special
Hermite semigroup and denoted by $e^{-tL}, t>0$) by the expansion
$$e^{-tL}f=(2\pi)^{-n}\sum_{k=0}^{\infty}e^{-(2k+n)t}f\times \varphi_k$$
for $f\in L^2(\mathbb{C}^n).$ Again, on a dense subspace, $e^{-tL}$ can be explicitly written as
$e^{-tL}f(z)=f\times p_t (z),$ where $$p_t(z)=(2\pi)^{-n}(\sinh t)^{-n}
e^{-\frac{1}{4}\coth t|z|^2}.$$ We shall identify $\mathbb{C}^n$ with $\mathbb{R}^{2n}$
and write $z=x+iy\in \mathbb{C}^n$ as $(x,u)\in \mathbb{R}^{2n}.$
For $f\in L^2(\mathbb{C}^n),$ the function $e^{-tL}f(x,u)$ is given by
$$f\times p_t(x,u)=\int_{\mathbb{R}^n} f(x',u')p_t(x-x',u-u')e^{-\frac{i}{2}(x'u-xu')}dx' du'.$$
The special Hermite functions $\Phi_{\alpha \beta}(x,u)$ can be extended to
$\mathbb{C}^{2n}$ as entire functions $\Phi_{\alpha \beta}(z,w)$. Moreover,
 $f\times p_t(x,u)$  also extends to $\mathbb{C}^{2n}$ as an entire function
$f\times p_t(z,w),$ where $z=x+iy, w=u+iv.$

Let $\mathcal{B}_t^*(\mathbb{C}^{2n})$ denote the set of all
entire functions on $\mathbb{C}^{2n}$ which are square integrable
with respect to the weight function
$$ W_t(z,w)=W_t(x+iy,u+iv)=4^n e^{(uy-vx)}p_{2t}(2y,2v).$$ Then the
following result is known, see \cite {Kro05}:
\begin{theorem}
The special Hermite semigroup
$e^{-tL}:L^2(\mathbb{R}^{2n})\rightarrow
\mathcal{B}_t^*(\mathbb{C}^{2n})$ is an isometric isomorphism.
\end{theorem}
 Let $\widetilde{\Phi}_{\alpha
\beta}(z,w)= e^{-(2|\beta|+n)t}\Phi_{\alpha \beta}(z,w).$ Then the family
$\{\widetilde{\Phi}_{\alpha \beta}/ \alpha, \beta \in
\mathbb{N}^n\}$ forms an orthonormal basis for
$\mathcal{B}_t^*(\mathbb{C}^{2n}).$ For further details, we refer to
\cite{Tha06}.

Laguerre Sobolev spaces were introduced by Peetre and Sparr in 1975.
These Sobolev spaces along with Hermite Sobolev spaces were studied
in connection with the regularity of the twisted spherical means in
\cite {Tha993}.

The Laguerre Sobolev space $W_L^{m,2}(\mathbb{C}^n)$ is defined
to be the image of $L^2(\mathbb{C}^n)$ under $L^{-m},$ where $L$
denotes the special Hermite operator mentioned earlier. In other
words $f\in W_L^{m,2}(\mathbb{C}^n)$ if and only if
$$||f||_{W_L^{m,2}}^2=\sum_\alpha\sum_{\beta}(2|\beta|+n)^{2m}|(f,\Phi_{\alpha \beta})_{L^2(\mathbb{C}^n)}|^2 <\infty.$$
The holomorphic Sobolev space
$W_t^{*~m,2}(\mathbb{C}^{2n})$ is defined as the image of
$W_L^{m,2}(\mathbb{C}^{n})$ under $e^{-tL}.$ Recall that
$\widetilde{\Phi}_{\alpha \beta}=e^{-tL}\Phi_{\alpha \beta}.$ The
norm in $W_t^{*~m,2}(\mathbb{C}^{2n})$ is given by
$$||F||_{W_t^{*~m,2}(\mathbb{C}^{2n})}=\sum_\alpha\sum_{\beta}(2|\beta|+n)^{2m}|(F,\tilde{\Phi}_{\alpha \beta})_{B_t^*(\mathbb{C}^{2n})}|^2.$$
Let $G_t^m(\mathbb{C}^{2n})$ denote the space of all  holomorphic
functions on $\mathbb{C}^{2n}$ which are square integrable with
respect to $|\frac{d^{2m}}{dt^{2m}}W_t(z,w)| dz dw .$ As in
section 2, we can show that
\begin{eqnarray}
(F,G)=\int_{\mathbb{C}^{2n}}F(z,w)\overline{G(z,w)}\frac{d^{2m}}{dt^{2m}}W_t(z,w)dzdw
\end{eqnarray}
is an inner product on $G_t^m(\mathbb{C}^{2n}).$ Let
$\mathcal{B}_t^{*,m}(\mathbb{C}^{2n})$ be the completion of
$G_t^m(\mathbb{C}^{2n})$ with this inner product (3.1). We have the following
result.

\begin{theorem}
For every non negative integer $m, W_t^{*~m,2}(\mathbb{C}^{2n})$
coincides with $\mathcal{B}_t^{*,m}(\mathbb{C}^{2n})$ and the special Hermite
semigroup $e^{-tL}$ is an isometric isomorphism of
$W_L^{m,2}(\mathbb{C}^n)~\mbox{onto}~                                           \mathcal{B}_t^{*,m}(\mathbb{C}^{2n})$ up to a constant
multiple.
\end{theorem}

This theorem is proved by  using
the fact that $F$ can be expressed in two ways $F=\sum (F,\tilde
\Phi_{\alpha \beta})_{B_t^*(\mathbb{C}^{2n})}\tilde \Phi_{\alpha \beta}, ~
F(z,w)=\sum_{\alpha,\beta}(F,\Phi_{\alpha \beta})_{L^2(\mathbb{C}^{n})}\Phi_{\alpha
\beta}$ where $z=x+iy, w=u+iv$ and $$(F,\Phi_{\alpha
\beta})_{L^2(\mathbb{C}^n)}=\int_{\mathbb{C}^n}F(x+iy)\overline{\Phi_{\alpha \beta}(x+iy)}dxdy,$$
the orthogonality relation
$$\int_{\mathbb{C}^{2n}}\Phi_{\alpha \beta}(z,w)\overline{\Phi_{\mu \nu}(z,w)} W_t(z,w)dzdw
=e^{2(2|\beta|+n\nu)t}\delta_{\alpha \mu}\delta_{\beta \nu},$$
and the fact that
$\widetilde{\Phi}_{\alpha \beta}\in
\mathcal{B}_t^{*,m}(\mathbb{C}^{2n}),$ and proceeding as in
section 2.

\section{The Image of Schwartz class functions under Hermite and special Hermite semigroups}

First, we shall describe the image of $\mathcal{S}(\mathbb{R}^n)$
under $e^{-tH}.$ In order to do this, first we shall
obtain pointwise estimates for a function $F\in H_t(\mathbb{C}^n);$
i.e for $ F = e^{-tH}f,~ f\in L^2(\mathbb{R}^n).$ Since $F\mapsto F(z)$
is a continuous linear functional on $H_t(\mathbb{C}^n)$ for each
$z\in \mathbb{C}^n,$ an application of Riesz representation theorem shows
that there exists a unique $K_t(z,.)\in H_t(\mathbb{C}^n)$ such that
\begin{eqnarray*}
F(z)&=&(F, K_t(z,.))\\&=&
\int_{\mathbb{C}^n}F(w)\overline{K_t(z,w)}U_t(w)dw.
\end{eqnarray*}
The function  $ K_t(z,w)$ is called the reproducing kernel
for $H_t(\mathbb{C}^n).$ By expanding $F$ in terms of
$\widetilde{\Phi}_{\alpha},$ we can write
\begin{eqnarray}
F(z)&=&\nonumber \int_{\mathbb{C}^n}
F(w)\sum_{\alpha}{\overline{\widetilde{\Phi}_{\alpha}(w)}}\widetilde{\Phi}_{\alpha}(z)U_t(w)dw\\&=&
\int_{\mathbb{C}^n} F(w)\overline{K_t(z,w)}U_t(w)dw.
\end{eqnarray}
This means that
$$K_t(z,w)= \sum_{\alpha}e^{-(2|\alpha|+n)2t}\Phi_{\alpha}(w)
\Phi_{\alpha}(\overline{z}).$$
Cauchy-Schwarz inequality applied to (4.1) gives us
\begin{eqnarray}
|F(z)|^2\leq
||F||^2  ||K_t(z,.)||^2 = ||F||^2 K_t(z,z).
\end{eqnarray}
By using Mehler's formula, we can explicitly calculate $ K_t(z,z)$. In fact,
\begin{eqnarray}
K_t(z,z)=(2\pi)^{-\frac{n}{2}}(\sinh4t)^{-\frac{n}{2}}
e^{-\frac{1}{2}\coth4t(z^2+\overline{z}^2)}e^{\frac{1}{\sinh4t}z\overline{z}}
\end{eqnarray}
Since
$$-\coth4t(x^2-y^2)+\mbox{cosech}4t(x^2+y^2)=-x^2\tanh2t+y^2\coth2t$$
we obtain
\begin{eqnarray}
|F(z)|^2\leq
 C_n (\sinh4t)^{-\frac{n}{2}}e^{-\tanh2t {x^2}+ \coth2t~{y^2}} ||F||^2.
\end{eqnarray}
Thus we have obtained a pointwise estimate for functions
$F\in H_t(\mathbb{C}^n).$

In order to obtain  pointwise estimates for $F\in
W_t^{m,2}(\mathbb{C}^n)$ we observe that the reproducing kernel for $
W_t^{m,2}(\mathbb{C}^n)$ is given by
$$ K_t^{2m}(z,w) = \sum_{\alpha}(2|\alpha|+n)^{-2m}\overline{
\widetilde{\Phi}_{\alpha}(z)} \widetilde{\Phi}_{\alpha}(w).$$ We can write
this as
$$ K_t^{2m}(z,w) = \frac{1}{(2m-1)!}\int_0^{\infty}s^{2m-1}K_{s+t}(z,w) ds.$$
Using the explicit formula for $ K_s(z,z) $ we have
$$
K_t^{2m}(z,z) = \frac{(2\pi)^{-\frac{n}{2}}}{(2m-1)!}\int_0^{\infty}s^{2m-1}(\sinh4(t+s))^{-\frac{n}{2}}$$
$$ \times e^{-x^2(\coth4(t+s)-\mbox{\tiny cosech}4(t+s))}
e^{y^2(\mbox{\tiny cosech} 4(t+s)+\coth4(t+s))}ds$$
$$ = \frac{(2\pi)^{-\frac{n}{2}}}{(2m-1)!}\int_0^{\infty}s^{2m-1}(\sinh4(t+s))^{-\frac{n}{2}}
e^{-x^2\tanh2(t+s)+y^2(\coth2(t+s))}ds.$$

From the above expression for the reproducing kernel it is now an easy matter
to establish the following pointwise estimates for functions from the
holomorphic Sobolev spaces.

\begin{theorem}(Sobolev-embedding theorem)
Let $ m $ be a non-negative integer. Then every $ F
\in W_t^{m,2}(\mathbb{C}^n)$ satisfies the estimate
$$|F(z)|^2\leq C(1+x^2+y^2)^{-2m}e^{-x^2\tanh 2t+y^2 \coth 2t}.$$
\end{theorem}
\begin{proof} In order to prove the theorem we need to estimate the integral
appearing in the representation of the reproducing kernel. We rewrite the
kernel  as $ C_t  e^{-x^2\tanh 2t+y^2\coth 2t} I $ where
$$
I =  \int_0^{\infty}s^{2m-1}
(\sinh4(t+s))^{-\frac{n}{2}} e^{-x^2 (\tanh2(t+s)-\tanh 2t)}$$
$$ \times  e^{y^2 (\coth2(t+s)-\coth 2t)}ds $$ which after some
simplification yields
\begin{eqnarray*}
I &=& \int_0^{\infty}s^{2m-1}(\sinh4(t+s))^{-\frac{n}{2}}
e^{- y^2 \frac{\sinh 2s}{\sinh 2(t+s)\sinh(2t)}- x^2\frac{\sinh 2s}
{\cosh 2(t+s)\cosh(2t)}}ds.
\end{eqnarray*}
Thus we only need to show that the above integral is bounded by a constant
times $ (1+x^2+y^2)^{-2m}.$

To prove this estimate we break up the above integral into two parts. Using
the fact that $\sinh $ and $ \cosh $ are increasing functions and
$ \sinh s > s $ we see that
$$
\int_0^{t}s^{2m-1}(\sinh4(t+s))^{-\frac{n}{2}}
e^{- y^2 \frac{\sinh 2s}{\sinh 2(t+s)\sinh(2t)}- x^2 \frac{\sinh 2s}
{\cosh 2(t+s)\cosh(2t)}}ds
$$
is bounded by
$$ \int_0^{\infty} s^{2m-1} e^{-2ns} e^{-2(\frac{x^2}{\cosh^2 4t}+\frac{y^2}
{\sinh^2 4t})s} ds \leq C_t (1+x^2+y^2)^{-2m}.$$
On the other hand the integral
$$
\int_t^{\infty}s^{2m-1}(\sinh4(t+s))^{-\frac{n}{2}}
e^{- y^2 \frac{\sinh 2s}{\sinh 2(t+s)\sinh(2t)}- x^2 \frac{\sinh 2s}
{\cosh 2(t+s)\cosh(2t)}}ds
$$
is bounded by a constant times
$$ e^{- (a_t x^2 +b_t y^2)} \int_0^{\infty}s^{2m-1} e^{-2ns} ds $$ where
$ a_t $ and $ b_t $ are the infima of $ \frac{\sinh 2s}{\cosh 2(t+s)\cosh(2t)}
$ and $ \frac{\sinh 2s}{\sinh 2(t+s)\sinh(2t)} $ over $ s > t $ respectively.
The above clearly gives the required estimate.
\end{proof}

Now, we are in a position to prove the following result which characterises
the image of $\mathcal{S}(\mathbb{R}^n)$ under $e^{-tH}.$
\begin{theorem}
Let $ t> 0 $ be fixed. Suppose $F$ is a holomorphic function on
$\mathbb{C}^n.$ Then there exists a function $f\in
\mathcal{S}(\mathbb{R}^n)$ such that  $F=e^{-tH}f$ iff $F$ satisfies
\begin{equation} |F(z)|^2\leq A_m \frac{e^{-(\tanh 2t) x^2+ (\coth 2t)
y^2}}{(1+x^2+y^2)^{2m}}
\end{equation}
for some constants $A_m, m=1,2,3,\cdots.$
\end{theorem}
\begin{proof}
If $f\in \mathcal{S}(\mathbb{R}^n),$ then $f\in
W_{H}^{m,2}(\mathbb{R}^n)~\forall ~m,$ which in turn implies that
$F=e^{-tH}\in W_t^{m,2}(\mathbb{C}^n)~\forall ~m.$ Then $(4.5)$
follows from theorem 4.1. Conversely, suppose $F$ satisfies
$(4.5).$ Then by choosing $ m $ large we see that
\begin{equation}
\int_{\mathbb{C}^n}|F(z)|^2 U_t(z)dz <\infty
\end{equation}
from which it follows that $F\in H_t(\mathbb{C}^n).$ Thus there
exists a function $f\in L^2(\mathbb{R}^n)$ such that $F=e^{-tH}f.$
Since $\frac{d^{2m}}{dt^{2m}}U_t(z)$ is the sum of $2m+1$ terms,
where each term is of the form $(p(t)x^2+q(t)y^2+c)^k\leq C_t
(1+x^2+y^2)^{2m},$ with $k\leq 2m.$ Thus it follows from (4.5)
that $F\in
\mathcal{B}_t^{m}(\mathbb{C}^n)=W_t^{m,2}(\mathbb{C}^n)$ using
Theorem 2.3. This leads to the fact that $F\in
W_t^{m,2}(\mathbb{C}^n) ~\forall ~m.$ Consequently $f\in
W_{H}^{m,2}(\mathbb{R}^n)~ \forall~ m.$ But since $\bigcap_m
W_{H}^{m,2}(\mathbb{R}^n)=\mathcal{S}(\mathbb{R}^n),$ the result
follows.
\end{proof}

Now we shall characterise the image of ${\mathcal{S}(\mathbb{C}^n)}$ under $e^{-tL}.$
In order to do so, first we get  pointwise estimates for functions in
$\mathcal{B}_t^*(\mathbb{C}^{2n}).$
Let $F\in \mathcal{B}_t^*(\mathbb{C}^{2n}).$ Then
$|F(z,w)|=|e^{-tL}f(z,w)|=|(f\times p_t)(z,w)|$
where
\begin{eqnarray*}
p_t(x,u) =  (2\pi)^{-n} (\sinh t)^{-n}e^{-\frac{1}{4}\coth t(x^2+y^2)}
\end{eqnarray*}
(see \cite {Tha06} for details.) Recalling the definition of twisted convolution,
$$ f\times p_t
(x,u)=\int_{\mathbb{R}^{2n}}f(x',u')p_t(x-x',u-u')e^{-\frac{i}{2}(x'u-xu')}
dx'du',$$ we get
\begin{eqnarray*}
|f\times p_t
(z,w)|&=&|\int_{\mathbb{R}^{2n}}f(x',u')p_t(z-x',w-u')e^{-\frac{i}{2}(x'w-zu')}dx'du'|
\\&\leq& \int_{\mathbb{R}^{2n}}|f(x',u')||p_t(z-x',w-u')|
e^{\frac{x'v}{2}+\frac{y'u}{2}}dx'du'
\\&\leq&||f||_{L^2(\mathbb{C}^n)}(2\pi)^{-n}(\sinh{2t})^{-n}e^{vx-uy}
e^{\coth{4t(y^2+v^2)}}
\end{eqnarray*}
where $w=u+iv, z=x+iy.$
Proceeding as in section 3, we can show that if $F\in
W_t^{*~m,2}(\mathbb{C}^{2n}),$ then
\begin{eqnarray}\label{equ4.6}|F(z,w)|^2\leq C_t\frac{ e^{vx-uy}e^{\coth{4t}(y^2+v^2)}}{(1+y^2+v^2)^{2m}}.
\end{eqnarray}
But we can do better than this for Schwartz functions. We require the
following simple lemma.

%**** I do not understand the following arguments****
\begin{lemma}
If $f\in \mathcal{S}(\mathbb{R}^{2n}),$ then
\begin{eqnarray}
e^{-tL}(\prod   _{j=1}^n (\frac{\partial}{\partial x_j}-ax_j)f) = 
\prod   _{j=1}^n (-az_j+bw_j)e^{-tL}f  \\
e^{-tL}( \prod   _{j=1}^n (\frac{\partial}{\partial u_j}+bu_j)f ) =
\prod   _{j=1}^n (bz_j+aw_j)e^{-tL}f
\end{eqnarray}
where $a=-\frac{1}{2}\coth t, ~ b=\frac{i}{2}.$
\end{lemma}
\begin{proof} We shall prove the result for $n=1.$
Consider
$$
e^{-tL}(\frac{\partial}{\partial x'}f)(z,w) = $$
$$ \frac{1}{4\pi}(\sinh t)^{-1}\int_{\mathbb{R}^2}
(\frac{\partial}{\partial x'}f)(x',u')e^{-\frac{1}{4}\coth t[(z-x')^2+(w-u')^2]}e^{-\frac{i}{2}(x'w-zu')}dx' du'.$$
Integration by parts leads to
$$
e^{-tL}(\frac{\partial}{\partial x'}f)(z,w) $$
$$ = -\frac{1}{4\pi}(\sinh t)^{-1}\int_{\mathbb{R}^2}
f(x',u')\frac{\partial}{\partial x'}[e^{-\frac{1}{4}\coth t[(z-x')^2+(w-u')^2]}e^{-\frac{i}{2}(x'w-zu')}]dx' du'$$
$$ =  -\frac{1}{4\pi}(\sinh t)^{-1}\int_{\mathbb{R}^2}
f(x',u')[\frac{1}{2}\coth t(z-x')-\frac{i}{2}w]$$
$$ \times e^{-\frac{1}{4}\coth t[(z-x')^2+(w-u')^2]} e^{-\frac{i}{2}(x'w-zu')}dx' du'$$
$$ = (-\frac{1}{2}z\coth t +\frac{i}{2}w)e^{-tL}f(z,w)
+\frac{1}{8\pi}(\sinh t)^{-1}\coth t $$
$$ \times \int_{\mathbb{R}^2}x'f(x',u')e^{-\frac{1}{4}\coth t[(z-x')^2+(w-u')^2]}e^{-\frac{i}{2}(x'w-zu')}dx' du'.$$
Thus
$$
e^{-tL}(\frac{\partial}{\partial x'}f)(z,w) = $$
$$ (-\frac{1}{2}z \coth t +\frac{i}{2}w)e^{-tL}f(z,w)+\frac{1}{2}\coth t ~e^{-tL}(x'f)(z,w).$$
Hence
$$e^{-tL}\frac{\partial}{\partial x'}f=(-az+bw)e^{-tL}f+a e^{-tL}(x'f).$$
So $$e^{-tL}(\frac{\partial}{\partial x'}-ax')f=(-az+bw)e^{-tL}f.$$
In the case of $n$ dimension if $x'=(x_1',x_2',\cdots,x_n'),$ then we have to apply the same procedure
for each $x_j', j=1,2,\cdots, n $ in order to obtain $(4.8).$ In a similar 
way $(4.9)$ can be proved.
\end{proof}

In what follows we use the standard vector notation 
$ (\frac{\partial}{\partial x'}- ax')^\alpha =\prod_{j=1}^n 
(\frac{\partial}{\partial x_j'}- ax_j')^{\alpha_j} $ where $ \alpha $ is a 
multi-index.

 \begin{theorem}
 Suppose $F$ is a holomorphic function on $\mathbb{C}^{2n}.$ Fix
 $t>0.$ Then there exists a function
$f\in \mathcal{S}(\mathbb{R}^{2n})~\mbox{with}~F=e^{-tL}f$ iff $F$ satisfies
$$|F(z,w)|^2\leq B_m\frac {e^{vx-uy}e^{\coth{4t}(y^2+v^2)}}
{(1+x^2+y^2+u^2+v^2)^{2m}}$$
 for some constants $B_m, m=1,2,3,\cdots.$
 \end{theorem}
\begin{proof}
Multiplying $(4.8)$ by $a$, multiplying $(4.9)$ by $b$ and subtracting $(4.8)$ from $(4.9)$ we get
\begin{eqnarray}
&&-a e^{-tL}[(\frac{\partial}{\partial x}-ax)f]+ b e^{-tL}[(\frac{\partial}{\partial u}+bu)f] =(a^2+b^2)z e^{-tL}f.
\end{eqnarray}
Multiplying $(4.8)$ by $b$,  multiplying $(4.9)$ by $a$ and adding $(4.8)$ and  $(4.9)$ we get
\begin{eqnarray}
&&b e^{-tL}[(\frac{\partial}{\partial x}-ax)f]+ a e^{-tL}[(\frac{\partial}{\partial u}-bu)f] =(a^2+b^2)w e^{-tL}f.
\end{eqnarray}
Applying the lemma $(4.3)$ iteratively we obtain the following :
Let $k=(k_1,k_2,...,k_n), l=(l_1,l_2,...,l_n)\in \mathbb{N}^n$ and $x^k={x_1}^{k_1}{x_2}^{k_2}...{x_n}^{k_n}.$
If $f\in \mathcal{S}(\mathbb{R}^{2n})$ then $(\frac{\partial}{\partial x}-ax)^k f,(\frac{\partial}{\partial u}+bu)^l f
\in \mathcal{S}(\mathbb{R}^{2n})$ and $$e^{-tL}[(\frac{\partial}{\partial x}-ax)^k f]=(-az+bw)^k e^{-tL} f$$
$$e^{-tL}[(\frac{\partial}{\partial u}+bu)^l f]=(bz+aw)^l e^{-tL} f.$$
As in the case of (4.11) and (4.12), by carrying out appropriate algebraic manipulation we can find a differential operator $T_{k,l}$ such that $$e^{-tL}(T_{k,l}f)=z^kw^l e^{-tL}f$$ for any $k,l\in \mathbb{N}^n.$ Since $T_{k,l}f\in \mathcal{S}(\mathbb{R}^{2n}),$ it follows from (\ref{equ4.6}) that
\begin{eqnarray*}
|(z^kz^l e^{-tL}f)(z,w)|^2 &\leq&
C_t ~e^{vx-uy}e^{\coth4t(y^2+v^2)}
\end{eqnarray*}
Thus $$|(1+|z|^2+|w|^2)^m F(z,w)|^2\leq
C_{t,m}e^{vx-uy}e^{\coth4t(y^2+v^2)},$$ where $F=e^{-tL}f.$
\end{proof}

\section{ Tempered distributions and a Paley-Wiener theorem for the windowed
Fourier transform}

In this section we consider the image of tempered distributions on
$ \mathbb{R}^n $ under the Hermite semigroup. The characterisation obtained
leads to a
Paley-Wiener theorem for the windowed Fourier transform. We prove the
following analogue of Theorem 4.2 for tempered distributions.

\begin{theorem}
 Suppose $F$ is a holomorphic function on $\mathbb{C}^{n}.$ Then there
exists a distribution $f\in \mathcal{S}'(\mathbb{R}^{n})~\mbox{with}~F=e^{-tH}f$
if and only if $F$ satisfies
$$|F(z)|^2\leq C (1+|z|^2)^{2m} e^{-x^2 \tanh 2t + y^2 \coth 2t} $$
for some non-negative integer $ m.$
\end{theorem}

We obtained  Theorem 4.2 as a consequence of Theorem 4.1 using the fact that
the intersection of all the Hermite Sobolev spaces $ W^{m,2}_H(\mathbb{R}^n)$
is precisely the space of Schwartz functions. Since the union of all
$ W^{m,2}_H(\mathbb{R}^n)$ is $ \mathcal{S}'(\mathbb{R}^{n}) $ we only need
to prove the following analogue of Theorem 4.1 for functions from
$ W^{-m,2}_t(\mathbb{C}^n)$ where $ m $ is a positive integer.

\begin{theorem} Let $ m $ be a positive integer. Then every
$ F \in W^{-m,2}_t(\mathbb{C}^n) $ satisfies the estimate $$
|F(z)|^2 \leq C (1+|z|^2)^{2m}e^{-x^2 \tanh 2t + y^2 \coth 2t}.$$ Conversely,
if an entire function $ F $ satisfies the above estimate, then $ F $ belongs
to $ W^{-m-n,2}_t(\mathbb{C}^n) $
\end{theorem}

The necessity of the pointwise estimate is easy to establish. In fact, the
reproducing kernel for $ W^{-m,2}_t(\mathbb{C}^n) $ is given by
$$ K_t^{-2m}(z,w) = \sum_{\alpha} (2|\alpha|+n)^{2m}e^{-2(2|\alpha|+n)t}
\Phi_\alpha(\bar{z})\Phi_\alpha(w).$$ Therefore, we only need to estimate
the $ (2m)$-th derivative of $ K_t(z,z) $ with respect to $ t.$ It is
easy to see that this leads to the estimate
$$ |F(z)|^2 \leq C (1+|z|^2)^{2m}e^{-x^2 \tanh 2t + y^2 \coth 2t}.$$

To prove the converse, we need to make use of duality. As
$ W^{m,2}_H(\mathbb{R}^n)$ and $ W^{-m,2}_H(\mathbb{R}^n)$ are dual to each
other, it follows from the definition that $ W^{m,2}_t(\mathbb{C}^n) $ and
$ W^{-m,2}_t(\mathbb{C}^n) $ are dual to each other. The duality bracket is
given by
$$ (F,G) = \int_{\mathbb{C}^n} F(z) \overline{G(z)} U_t(z) dz.$$ We refer to
\cite {Tha07} for details of the duality argument in the context of compact
symmetric spaces. As the same proof works in our situation as well we do not
give any details but only a brief sketch. If $ F $ satisfies the given estimates then for any $ G \in
 W^{m+n,2}_t(\mathbb{C}^n) $ the integral defining $ (F,G) $ converges and
hence $ F $  defines a continuous linear functional on
$ W^{m+n,2}_t(\mathbb{C}^n).$ Consequently, $ F $ belongs to
$ W^{-m-n,2}_t(\mathbb{C}^n) $ which proves the converse.

We remark that the holomorphic Sobolev spaces $ W^{-m,2}_t(\mathbb{C}^n) $
are weighted Bergman spaces with a nonnegative weight function which can be
given a representation in terms of Riemann-Liouville fractional integrals,
see \cite {Tha07} for details.

We now show that the characterisation obtained in Theorem 5.1 leads to a
Paley-Wiener theorem for the windowed Fourier transform. This transform, also
known as short-time Fourier transform is used very much in Gabor analysis.
Given a Schwartz class function $ g $ the windowed Fourier transform of
$ f $ (with window $ g $) is defined by
$$ V_g(f)(x,y) = (2\pi)^{-\frac{n}{2}} \int_{\mathbb{R}^n} f(u)g(u-y)
e^{-i x \cdot u} du.$$ See \cite {Gro01} and \cite{Gro04} for the study of
windowed Fourier transforms in connection with Hardy's theorem and
characterisation of Schwartz functions. We are mainly interested in the case
where $ g(x) = \varphi_a(x) = c_n e^{-\frac{1}{2}a |x|^2}.$ More precisely, we consider
the transform
$$ T_a(f)(x) = V_{\varphi_a}(f)(x,0).$$ Note that  $ T_a(f) $ is
well defined as a function even if $ f $ is a tempered distribution.

The classical Paley-Wiener theorem characterises the space of all compactly
supported distributions in terms of the holomorphic properties of their
Fourier transforms. Even though we can define Fourier transforms of tempered
distributions we cannot hope for any such characterisation since the space of
all tempered distributions is invariant under the Fourier transform. On the
other hand we see that $ T_a f $ extends to $ \mathbb{C}^n $ as an entire
function even when $ f $ is only a tempered distribution. This property of the
windowed Fourier transform allows us to prove the following analogue of
Paley-Wiener theorem.

\begin{theorem} For any $ a > 0 $ the windowed Fourier transform $ T_af(x) $
of a tempered distribution $ f $ on $ \mathbb{R}^n $ extends to
$ \mathbb{C}^n $ as an entire function which satisfies the estimate
$$ |T_af(x+iy)| \leq C (1+x^2+y^2)^m e^{\frac{1}{2}a^{-1}y^2} $$ for some non-
negative integer $ m.$ Conversely, if an entire function $ F $ satisfies such
an estimate, then $ F = T_af $ for a tempered distribution $ f $.
\end{theorem}
\begin{proof} This theorem is easily proved by relating the windowed Fourier
transform $ T_af $ with $ e^{-tH}f .$  Indeed,
considering the case $ a > 1 $ first and writing $ a = \coth(2t) $ for some
$ t > 0 $ we can easily verify that
$$ e^{-tH}f(z) = e^{-\frac{1}{2}\coth(2t)z^2}T_af(\frac{iz}{\sinh(2t)}) $$
for all
$ z \in \mathbb{C}^n.$ We obtain the required estimate on $ T_af(z) $ by
appealing to Theorem 5.1. Conversely, if $ F $ satisfies the given estimates
then again by Theorem 5.1 the function
$$ G(z) = e^{-\frac{1}{2}\coth(2t)z^2}F(\frac{iz}{\sinh(2t)}) $$
should be of the form $ e^{-tH}f(z) $ for a tempered distribution $ f.$

When $ a < 1 $ we take $ t > 0 $ so that $ a = \tanh(2t).$ Now the proof
requires an analogue of Theorem 5.1 for functions of the form
$ e^{-(t+i\frac{\pi}{4})H}f.$ But the image of tempered distributions under
$ e^{-(t+i\frac{\pi}{4})H} $ can be characterised in a similar way. The final
estimates do not depend on the factor $ e^{-i\frac{\pi}{4}H} $ which is just
the Fourier transform. This completes the proof of the theorem.
\end{proof}

We remark that the case $ a =1 $ can be read out from the results of \cite {
Bar67}. This case can be obtained as a limiting case of our results and we
think the proof presented here is simpler than the one given in \cite {Bar67}.
Finally, we remark that  we also have the following result which characteries
the image of compactly supported distributions under the Hermite semigroup.

\begin{theorem} Let $ f $ be a distribution supported in a ball of radius
$ R $ centered at the origin. Then for any $ t > 0 $ the function $ e^{-tH}f $
extends to $ \mathbb{C}^n $ as an entire function which satisfies
$$ |e^{-tH}f(z)| \leq C e^{-\frac{1}{2}\coth(2t)(x^2-y^2)}e^{\frac{R|x|}
{\sinh(2t)}}.$$
Conversely, any entire function $ F $ satisfying the above will be of the form
$ e^{-tH}f $ where $ f $ is supported inside a ball of radius $ R $ centered
at the origin.
\end{theorem}

The proof of this follows from the classical Paley-Wiener theorem for
compactly  supported distributions. We leave the details to the reader.

\bibliographystyle{amsplain}

\end{document}